\documentclass[leqno,12pt]{article} 
\usepackage{amsmath, amssymb}
\usepackage{amsthm} 
\usepackage{braket}
\usepackage{mathtools}
\DeclarePairedDelimiter{\abs}{\lvert}{\rvert}
\usepackage{latexsym}
\def\qed{\hfill $\Box$}
\usepackage{cases}

\makeatletter
    
    \@addtoreset{equation}{section}
\makeatother

\theoremstyle{plain} 
\newtheorem{theorem}{\indent\sc Theorem}[section]

\theoremstyle{definition} 
\newtheorem{definition}[theorem]{\indent\sc Definition}
\newtheorem{remark}[theorem]{\indent\sc Remark}

\newcommand{\deldel}{\sqrt{-1}\partial \overline{\partial}}

\newcommand{\e}{\varepsilon}

\newcommand{\poly}{\bigtriangleup}

\makeatletter

\def\address#1#2{\begingroup
\noindent\parbox[t]{7.8cm}{
\small{\scshape\ignorespaces#1}\par\vskip1ex
\noindent\small{\itshape E-mail}
\/: #2\par\vskip4ex}\hfill
\endgroup}
\makeatother

\title{Remarks on modified Ding functional for toric Fano manifolds}

\author{
\textsc{Satoshi Nakamura$^{*}$}
} 

\date{}

\begin{document}

\maketitle


\footnote{ 
$^{*}$Partly supported by Grant-in-Aid for JSPS Fellowships for Young Scientists, Number 17J02783.
}

\begin{abstract}
We give a characterization of relative Ding stable toric Fano manifolds in terms of the behavior of the modified Ding functional.
We call the corresponding behavior of the modified Ding functional the pseudo-boundedness from below.
We also discuss the pseudo-boundedness of the Ding / Mabuchi functional of general Fano manifolds.  
\end{abstract}


\section{Introduction}
The existence problem of K\"ahler Einstein metrics for Fano manifolds was one of the central problems in K\"ahler Geometry.
The vanishing of the Futaki invariant is known as an obstruction to the existence of K\"ahler Einstein metrics.
Mabuchi \cite{M} introduced the notion of {\it Generalized K\"ahler Einstein metrics}, which is a generalization of K\"ahler Einstein metrics for Fano manifolds with non-vanishing Futaki invariant.

For toric Fano manifolds, a few criterions for the existence of Generalized K\"ahler Einstein metrics were established by
Yao \cite{Y} in terms of Geometric invariant theoretic stabilities, and by the author \cite{N} in terms of the properness of a functional on the space of K\"ahler metrics.
Very recently,  Li and Zhou \cite{LZ} generalized these criterion 
for Fano group compactifications which are generalizations of toric Fano manifolds.

In the following, we fix the notation for toric Fano manifolds to state Yao and the author's criterions, and to state the main theorem of the present article.
Let $X$ be an $n$-dimensional toric Fano manifold, and
$\poly$ be the open set  in $\mathbb{R}^n$ such that the closure $\overline{\poly}$ is the reflexive integral Delzant polytope corresponding to $X$.
Note that $0\in\mathbb{R}^n$ is the only integral point in $\poly$.
Let $(\mathbb{C}^*)^n = (S^1)^n \times \mathbb{R}^n$ be the open dense orbit in $X$,
and $\xi_i := \log |z_i|^2$ be the coordinate of $\mathbb{R}^n$, where $\{z_i\}$ is the standard coordinate of $(\mathbb{C}^*)^n$.

Let $\omega_0 \in 2\pi c_1(X)$ be an $(S^1)^n$-invariant reference K\"{a}hler metirc on $X$.
It is well-known that there exists a smooth convex function 
$\phi_0 = \phi_0(\xi_1, \dots, \xi_n)$ on $\mathbb{R}^n$ such that
$\omega_0 = \deldel \phi_0$ holds on $(\mathbb{C}^*)^n$.
Let $u_0$ be the Legendre dual of $\phi_0$, that is, 
\begin{equation*}
u_0(x)=\sup_{\xi\in\mathbb{R}^n}(x\cdot\xi-\phi_0(\xi)).
\end{equation*}
Then $u_0$ is a smooth convex function on $\poly$.
Let
\[
\mathcal{C}:= \Set{u \in C^0(\overline{\poly}) | 
u \text{ is convex on } \poly \text{ and, } 
u-u_0 \in C^{\infty}(\overline{\poly})
},
\]
and let
\[
\tilde{\mathcal{C}}:=\set{ u\in\mathcal{C} | u\geq u(0) =0}.
\]
Abreu \cite{Ab} showed that there is the bijection between $\mathcal{C}$ and the space of $(S^1)^n$-invariant K\"{a}hler metrics in $[\omega_0]=2\pi c_1(X)$.

For $u\in\mathcal{C}$, we define the {\it modified Ding functional} \cite{Y} $\mathcal{D}$ of $X$ by
\[
\mathcal{D}(u) = -\log\int_{\mathbb{R}^n} e^{-(\phi-\inf\phi)}d\xi - u(0) + \int_{\poly}u\cdot l dx,
\]
where $\phi$ is the Legendre dual of $u$ and $l$ is the unique affine linear function on $\poly$ such that
\[
- u(0) + \int_{\poly}u\cdot l dx = 0
\]
holds for any affine linear $u\in\mathcal{C}$.
The critical point of the modified Ding functional is the generalized K\"{a}hler Einstein metric \cite{Y}.
We also define the {\it relative Ding-Futaki invariant} \cite{Y} $\mathcal{I}$ of $X$ by
\[
\mathcal{I}(u)=- u(0) + \int_{\poly}u\cdot l dx.
\]

We define another invariant $\alpha_X$ of $X$ by 
\begin{equation}\label{Yao}
\alpha_X = \max_{\overline{\poly}}(1-\abs\poly l),
\end{equation}
where $\abs\poly$ is the volume $\int_{\poly}dx$.
\begin{remark}
Originally the invariant $\alpha_X$ was introduced by Mabuchi \cite{M} for general Fano manifolds as an obstruction to the existence of Generalized K\"ahler Einstein metrics. 
Mabuchi showed that if $X$ admits Generalized K\"ahler Einstein metrics then $\alpha_X <1$ holds.
For toric Fano manifolds, Yao \cite{Y} gave the above explicit formula \eqref{Yao}.
\end{remark}

The criterions by Yao and the auther for the existence of Generalized K\"ahler Einstein metrics for toric Fano manifolds are as follows:
\begin{theorem}\label{YN}{\sc (\cite{Y, N})}
Let $X$ be a toric Fano manifold.
The followings are equivalent.
\begin{enumerate}
\item $X$ admits a unique toric invariant generalized K\"{a}hler Einstein metic.
\item $\alpha_X<1$.
\item $X$ is uniform relative Ding stable \cite{Y, N}. 
Namely, there exists a constant $\lambda>0$ such that
\[
\mathcal{I}(u) \geq \lambda \int_{\poly} u dx
\]
holds for any $u \in \tilde{\mathcal{C}}$.
\item The modified Ding functional $\mathcal{D}$ of $X$ is proper. 
Namely, there exists an increasing function $\mu(r)$ on $\mathbb{R}$ with the property
$
\lim_{r\to\infty}\mu(r) = \infty
$
such that
\[
\mathcal{D}(u) \geq \mu \Bigl( \int_{\poly}udx \Bigr)
\]
holds for any $u\in\tilde{\mathcal{C}}$.
\end{enumerate}
\end{theorem} 

Yao \cite{Y} also introduced the notion of the relative Ding stability for toric Fano manifolds.
Therefore it is natural to ask how the functional $\mathcal{D}$ behaves in this case.
The following is our main theorem of this article.
\begin{theorem}\label{main}
Let $X$ be a toric Fano manifold.
The followings are equivalent.
\begin{enumerate}
\item \label{a} $\alpha_X \leq 1$.
\item \label{b} $X$ is relative Ding stable \cite{Y}. 
Namely,
\[
\mathcal{I}(u) \geq 0
\] 
holds for all $u\in\tilde{\mathcal{C}}$.
\item \label{c} For any $\e > 0$ there exists $C_{\e}>0$ such that
\[
\mathcal{D}(u) \geq - \e \int_{\poly}udx -C_{\e}
\]
holds for any $u \in \tilde{\mathcal{C}}$.
\end{enumerate}
\end{theorem}
Note that the equivalence between \eqref{a} and \eqref{b} was essentially given by Yao \cite{Y}.
The author's contribution is the discussion on the condition \eqref{c}.

The condition \eqref{c} in Theorem \ref{main} is weaker than the boundedness from below of $\mathcal{D}$ on $\tilde{\mathcal{C}}$.
In section \ref{quasi},
we generalize the condition \eqref{c} to the {\it pseudo-boundedness} from below (Definition \ref{quasi bounded})
for any functionals on the space of K\"ahler metrics of any K\"ahler manifolds.
Then we observe the following:
\begin{theorem}\label{main2}
Let $X$ be a Fano manifold.
Then the pseudo-boundedness from below of the Ding functional of $X$ implies the boundedness from below of itself.
The same statement holds for the Mabuchi functional.
\end{theorem}

{\sc Acknowledgments.}
The author would like to thank Professor Shigetoshi Bando and Doctor Ryosuke Takahashi for several helpful comments and constant encouragement.
He is partially supported by Grant-in-Aid for JSPS Fellowships for Young Scientists, Number 17J02783.

\section{Proof of Theorem \ref{main}.}
The proof is a slight modification of that of \cite[Theorem 1.1]{N}.

{\it Proof of} \eqref{a} $\Rightarrow$ \eqref{b}:
By Yao's formula \eqref{Yao}, for any $u\in\tilde{\mathcal{C}}$, we have
\[
\mathcal{I}(u) = \int_{\poly}u\cdot ldx \geq \frac{1-\alpha_X}{\abs\poly}\int_{\poly}udx \geq 0.
\]

{\it Proof of} \eqref{b} $\Rightarrow$ \eqref{c}:
For fixed $v_0 \in \mathcal{C}$ and its Legendre dual $\psi_0$, 
we define the bounded positive function $A$ on $\overline{\poly}$
by
\[
A(\nabla\psi_0) = \frac{e^{-\psi_0}}{\int_{\mathbb{R}^n}e^{-\psi_0}} \det(\nabla^2\psi_0)^{-1}.
\]
In fact $A$ is the exponential of the Ricci potential defined by $\psi_0$.
We also define the functional $\mathcal{D}_A$ by 
\[
\mathcal{D}_A(u)= -\log\int_{\mathbb{R}^n} e^{-(\phi-\inf\phi)}d\xi - u(0) + \int_{\poly}u\cdot A dx,
\]
where $\phi$ is the Legendre dual of $u$.
Note that the functional $\mathcal{D}_A$ can be defined on the space of bounded convex functions on $\overline{\poly}$, and is convex on this space \cite[Proposition 2.15]{BB}.
Note also that $v_0$ minimizes $\mathcal{D}_A$ on $\mathcal{C}$.
Indeed $v_0$ is the critical point of $\mathcal{D}_A$,
and $\mathcal{D}_A$ is convex on $\mathcal{C}$.
See \cite[Section 3.3]{Y} for more details.

Now we estimate the nonlinear term of $\mathcal{D}$.
Let us take any $\e\in (0,1]$, and any $u\in\tilde{\mathcal{C}}$. 
We denote the Legendre dual of $u$ by $\phi$.
Note that, by properties of the Legendre duality, $\inf\phi=-u(0) ( =0 )$ holds, 
and the Legendre dual of $\e u(x)$ is $\e\phi(\xi/\e)$.
Then we have
\begin{eqnarray}
-\log\int_{\mathbb{R}^n} e^{-(\phi-\inf\phi)}d\xi  
&\geq& -\log\int_{\mathbb{R}^n} e^{-\e\phi}d\xi  \label{e}\\ 
&=& -\log\int_{\mathbb{R}^n} e^{-\e\phi(\xi /\e)}d\xi + n\log \e \nonumber  \\
&=& \mathcal{D}_A(\e u) - \int_{\poly} \e u A dx + n\log \e \nonumber.
\end{eqnarray}
In order to estimate $\mathcal{D}_A(\e u)$, let us consider
\[
u':=\e^2u+(1-\e^2)u_0.
\]
Since $u'-u_0  \in C^{\infty}(\overline{\poly})$, the convex function $u'$ is in $\mathcal{C}$.
Thus $\mathcal{D}_A(u') \geq \mathcal{D}_A(v_0)$.
By the convexity of $\mathcal{D}_A$,
\[
\mathcal{D}_A(u') 
\leq \e \mathcal{D}_A(\e u) + (1-\e) \mathcal{D}_A((1+\e)u_0).
\]
It follows that
\begin{eqnarray}\mathcal{D}_A(\e u) 
&\geq& \frac{1}{\e}\mathcal{D}_A(u') - \frac{1-\e}{\e} \mathcal{D}_A((1+\e)u_0) \label{D} \\
&\geq& \frac{1}{\e}\mathcal{D}_A(v_0) - \frac{1-\e}{\e} \mathcal{D}_A((1+\e)u_0). \nonumber
\end{eqnarray}
By the assumption of the relative Ding stability, \eqref{e} and \eqref{D}, we have
\begin{eqnarray*}
\mathcal{D}(u) &=& -\log\int_{\mathbb{R}^n} e^{-(\phi-\inf\phi)}d\xi + \mathcal{I}(u) \\
&\geq& -\log\int_{\mathbb{R}^n} e^{-(\phi-\inf\phi)}d\xi \\
&\geq& 
 -\e\sup_{\poly}A\int_{\poly}udx +n\log\e 
 +\frac{1}{\e}\mathcal{D}_A(v_0) - \frac{1-\e}{\e} \mathcal{D}_A((1+\e)u_0).
\end{eqnarray*}
Replacing $\e\sup A$ by $\e$, we obtain the desired esitmate.

{\it Proof of} \eqref{c} $\Rightarrow$ \eqref{a}:
By Yao's formula \eqref{Yao} and the linearity of $l$, it suffice to show that 
$l(p)\geq 0$ for any vertex $p\in\partial{\poly}$.
As in the proof of \cite[Proposition 5.1]{N}, we can take a sequence of smooth convex function $\{v_i\}_i$ on $\overline{\poly}$ such that
(i) $v_i \geq v_i(0) = 0$, and
(ii) for fixed $K>0$, $v_i$ tends to the $K$ times of the Dirac function for the vertex $p\in\partial\poly$.
Then the convex function $u_i := \tilde{u}_0 + v_i$ is in  $\tilde{\mathcal{C}}$, 
where $\tilde{u}_0 \in \tilde{\mathcal{C}}$ is the normalization of
$u_0 \in \mathcal{C}$. 
Let 
\[
\phi_i(\xi):=\sup_{x\in\poly}(x\cdot\xi-u_i(x))
\]
be the Legendre dual of $u_i$.
Note that
\[
\inf_{\mathbb{R}^n}\phi_i=0 \quad\text{and}\quad \phi_i(\xi)\leq\sup_{x\in\poly}(x\cdot\xi),
\]
since $\inf_{\mathbb{R}^n}\phi_i = -u_i(0)$ by a property of the Legendre duality, and $u_i \geq u_i(0) = 0$ by the definition of $u_i$.
It follows that
\[
\log\int_{\mathbb{R}^n}e^{-(\phi_i-\inf\phi_i)} d\xi \geq \log\int_{\mathbb{R}^n}e^{-\sup_{x\in\poly}(x\cdot\xi)}d\xi.
\]  
By the condition \eqref{c},
we thus have
\begin{eqnarray*}
\int_{\poly}u_i\cdot l dx \geq -\e \int_{\poly}u_idx -C_{\e} 
 + \log\int_{\mathbb{R}^n}e^{-\sup_{x\in\poly}(x\cdot\xi)}d\xi.
\end{eqnarray*}
By taking $i\to\infty$, we have
\begin{eqnarray*}
K\cdot l(p)+\int_{\poly}l\cdot \tilde{u}_0 dx \geq - \e \Bigl(K+\int_{\poly}\tilde{u}_0 dx\Bigr) -C_{\e} 
+ \log\int_{\mathbb{R}^n}e^{-\sup_{x\in\poly}(x\cdot\xi)}d\xi.
\end{eqnarray*}
By taking $K$ large, we have
$
l(p) \geq -\e.
$
It follows that $l(p)\geq 0$, since $\e>0$ is arbitrary.
\qed
 \section{Pseudo-boundedness.}\label{quasi}
In this section, we introduce the notion of the quasi boundedness from below for functionals on the space of the K\"ahler metrics, 
and prove Theorem \ref{main2}.

Let $(X, \omega)$ be a $n$-dimensional compact K\"{a}hler manifold.
We denote its volume $\int_X \omega^n$ by $V$.
Let $G$ be any maximal compact subgroup of $\mathrm{Aut}(X)$.
If $\omega$ is $G$-invariant then we define the space of $G$-invariant K\"ahler metrics in $[\omega]$ by
\[
\mathcal{M}_G(\omega) 
=\Set{\phi\in C^{\infty}(X) | 
\omega_{\phi}:= \omega +\deldel\phi >0 \text{ and } 
\phi \text{ is } G \text{ -invariant}.
}.
\]
For any $\phi\in\mathcal{M}_G(\omega)$, we define the Aubin's I-functional $I$ and J-functional $J$ \cite{Au} of $X$ by
\begin{eqnarray*}
I(\phi) &=& \frac{1}{V} \int_X \phi (\omega^n - \omega_{\phi}^n), \\
J(\phi) &=& \frac{1}{V} \int_0^1 dt \int_X \dot{\phi_t}(\omega^n-\omega^n_{\phi_t}),
\end{eqnarray*}
where $\{\phi_t\}_{t\in[0,1]}$ is a path in $\mathcal{M}_G(\omega)$ connecting $0$ to $\phi$.
Functionals $I, I-J$ and $J$ are non-negative and, these are equivalent.
Namely,
\begin{equation}\label{IJ}
0 \leq I(\phi) \leq (n+1)(I(\phi)-J(\phi)) \leq nI(\phi).
\end{equation}

\begin{definition}\label{quasi bounded}
Let $f : \mathcal{M}_G(\omega) \to \mathbb{R}$ be a functional.
$f$ is {\it pseudo-bounded} from below if for any $\e>0$ there exists $C_{\e}>0$ such that
\[
f(\phi) \geq -\e J(\phi) -C_{\e}
\]
holds for any $\phi \in \mathcal{M}_G(\omega)$.
\end{definition}
Note that the condition of the pseudo-boundedness from below is weaker than that of boudedness from below, since the Aubin's J-functional is non-negative. 

\begin{remark}
For toric manifolds, the functional $u\mapsto \frac{1}{\abs{\poly}}\int_{\poly}udx$ on $\tilde{\mathcal{C}}$ is essentially same as the Aubin's J-functional \cite[Lemma 2.2]{ZZ}.
Therefore the condition \eqref{c} in Theorem \ref{main} is the pseudo-boundedness from below of the modified Ding functional.
\end{remark}

In the following, let $X$ be a $n$-dimensional Fano manifold, and $\omega\in2\pi c_1(X)$ be a $G$-invariant reference K\"ahler metirc. 
We denote its volume $\int_X \omega^n$ by $V$.
The Ricci potential $f_{\omega}$ for $\omega$ is the function satisfying
\[
\mathrm{Ric}(\omega)-\omega=\deldel f_{\omega} \quad\text{and}\quad \int_X e^{f_{\omega}}\omega^n = V.
\]
We define the Ding functional $F$ \cite{D} and the Mabuchi functional $M$ \cite{M1} on $ \mathcal{M}_G(\omega)$ by 
\begin{eqnarray*}
F(\phi) &=& - \frac{1}{V} \int_0^1 \dot{\phi_t} \omega_{\phi_t}^n - \log\Bigl(\frac{1}{V}\int_Xe^{f_{\omega}-\phi}\omega^n\Bigr), \\
M(\phi) &=& - \frac{1}{V}\int_0^1 dt \int_X \dot{\phi_t} (S(\omega_{\phi_t}) - n ) \omega_{\phi_t}^n,
\end{eqnarray*}
where $S(\omega)$ is the scalar curvature for $\omega$.

{\it Proof of the Theorem \ref{main2}.}
First it is easy to see that $M(\phi) \geq F(\phi) + \frac{1}{V}\int_X f_{\omega}\omega^n$.
Thus the pseudo-boundedness from below of $F$ implies that of $M$.
Then, for any $t\in[0,1)$, there exist $\delta, C >0$ such that 
\[
M(\phi)+(1-t)J(\phi) \geq \delta J(\phi) -C
\]
holds for any $\phi\in \mathcal{M}_G(\omega)$.
By \cite[Theorem 1]{S} (see also Remark \ref{Sz}), the greatest lower bound of the Ricci curvature $R(X)$ is equal to 1, that is, 
\begin{eqnarray*}
R(X)&:=&\sup_{t\in[0,1]}\{\exists \omega \in 2\pi c_1(X) \text{ such that } \mathrm{Ric}(\omega) > t\omega \} \\
       &=& 1.
\end{eqnarray*}
By \cite[Theorem 3]{L}, it follows that $F$ and $M$ are bounded from below.
\qed
\begin{remark}\label{Sz}
In \cite[Theorem 1]{S}, Sz\'{e}kelyhidi used the functional 
\[
J_{\eta}(\phi) = \frac{1}{V}\int_0^1 dt \int_X\dot{\phi_t}(\Lambda_{\omega_{\phi_t}}\eta - n)\omega_{\phi_t}^n \quad (\eta\in2\pi c_1(X)),
\]
instead of the J-functional.
However, as pointed out in \cite{S}, the functional $I-J$ is essentially same as $J_{\eta}$.  
By \eqref{IJ}, so is the functional $J$.
\end{remark}
\section{Concluding remark on the condition \eqref{c} in Theorem \ref{main}.}
In view of Theorem \ref{main2}, it is natural to hope that the pseudo-boundedness from below of the modified Ding functional for toric Fano manifolds implies the boundedness from below of itself.
Yao \cite{Y} and Nitta-Saito-Yotsutani \cite{NSY} suggest that a stronger statement may hold for toric Fano manifolds.  
In fact, by using the formula \eqref{Yao} with computers, Yao checked that $\alpha_X\leq 1$ implies $\alpha_X <1$ for any $2$-dimensional toric Fano manifolds.
Nitta, Saito and Yotsutani checked that the same statement holds for any $3$- and $4$-dimensional toric Fano manifolds.
Therefore, by Theorem \ref{YN} and Theorem \ref{main}, for toric Fano manifolds of dimension less than or equal to $4$, the pseudo-boundedness from below of the modified Ding functional implies the properness of itself.

\bigskip

\address{
Mathematical Institute \\ 
Tohoku University \\
Sendai 980-8578 \\
Japan
}
{satoshi.nakamura.r8@dc.tohoku.ac.jp}
\end{document}